\documentclass[
final
]{dmtcs-episciences}

\usepackage[utf8]{inputenc}
\usepackage{subfigure}

\usepackage[round]{natbib}
\usepackage{amsmath,amsthm, amssymb} 
\usepackage{graphicx}
\usepackage{hyperref}

\numberwithin{equation}{subsection}

\newtheorem{defn}{Definition}
\newtheorem{prob}[defn]{Problem}
\newtheorem{prop}[defn]{Proposition}
\newtheorem{thm}[defn]{Theorem}
\newtheorem{cor}[defn]{Corollary}
\newtheorem{lem}[defn]{Lemma}
\newtheorem{conj}[defn]{Conjecture}
\newtheorem{rmk}[defn]{Remark}
\newtheorem{clm}[defn]{Claim}

\newcommand{\R}{\mathbb{R}}

\newcommand{\remove}[1]{}
\newcommand{\gad}{${\cal G}_C(P,d)$}
\newcommand{\ga}{${\cal G}_C(P)$}
\newcommand{\size}[1]{\left| #1 \right|}

\newcommand{\doi}[1]{doi: \href{https://doi.org/#1}{\nolinkurl{#1}}}

\author[Arijit Bishnu et al.]{Arijit Bishnu\affiliationmark{1}
  \and Mathew Francis\affiliationmark{2}
  \and Pritam Majumder\affiliationmark{1}\thanks{I would like to thank the Indian Statistical Institute for the funding.}}
\title{Curves, points, incidences and covering}
\affiliation{
  Indian Statistical Institute, Kolkata, India\\
  Indian Statistical Institute, Chennai, India}
\begin{document}
\publicationdata
{vol. 28:3}
{2026}
{1}
{10.46298/dmtcs.16899}
{2025-11-11; 2025-11-11; 2026-06-19}
{2026-06-20}
\maketitle
\begin{abstract}
	\remove{Motivated by the concept of arboricity or linear arboricity in an undirected graph that seeks to find   structures, like non-crossing trees or paths, that covers edges of the graph, we pose certain problems related to points and curves in geometry.} Given a point set, mostly a grid in our case, we seek upper and lower bounds on the number of curves that are needed to \emph{cover} the point set. We say a curve \emph{covers} a point if the curve passes through the point. We consider such coverings by {monotonic curves, lines, orthoconvex curves, etc.} We also study a problem that is converse of the covering problem -- if a set of $n^2$ points in the plane is covered by $n$ lines then can we say something about the configuration of the points?
\end{abstract}

\section{Introduction}
\label{sec:intro}
\noindent
Let $C$ be a family of curves (e.g., circle, convex curves, etc.) and $P$ be a set of points in $\R^d$. A curve $c \in C$ \emph{covers} a point $p \in P$ if $p$ lies on the curve $c$. We say that $C$ covers $P$ if all points in $P$ are covered by the union $\bigcup C$ of all members of $C$. We will be interested in the minimum cardinality of $C$ that covers $P$. 

\medskip

To start with, let us concentrate on the following problem. Let $P$ be a set of points in $\mathbb{R}^2$ in general position and the goal is to figure out the number of simple curves needed to cover $P$. The solution is trivial -- sort the points based on their $x$-coordinates and join them from left to right; that is we need just one simple curve to cover $P$. As we move from a simple curve with no restrictions whatsoever, to a straight line, the problem becomes hard and deserves non-trivial solutions~\cite{KratschPR16,BoissonnatDGK18,AfshaniBDN16,LangermanM05}. This obviously gives rise to a natural question about what happens to this problem if we consider point sets with some special configuration, like grids vis-a-vis different kinds of simple curves like circles, convex curves, orthoconvex curves, etc. To bring the variety of different point sets and curves under a unifying framework, we propose the following definition of \emph{geometric covering number}.

\normalsize

\begin{defn}{(Geometric covering number)}
	The geometric covering number of a point set $P$ in $\R^d$ with respect to a curve type $C$ (like circle, convex curve, orthoconvex curve, etc.), denoted as ${\cal G}_C(P,d)$, is the minimum number of curves of type $C$ needed to \emph{cover} all the points in $P$. We say a curve \emph{covers} a point if the point lies on the curve. If the dimension is implicit, we will just write \ga~instead of \gad.
\end{defn}

\remove{
	The starting point of this work is an effort to port the notion of arboricity in graphs to geometry. The concept of arboricity in graphs~\cite{diestel} 
	is to cover the edges of the graph by some structures like tree or path. More specifically, \emph{arboricity} of an undirected graph is the minimum number of forests into which its edges can be partitioned. Alternately, it can be defined as the minimum number of spanning forests that can cover all the edges of the graph. Similarly, the \emph{linear arboricity}~\cite{A1,A2,Harary,Alon-arboricity} of a graph is defined as the minimum number of linear forests (i.e., a collection of paths) into which the edges of the graph can be partitioned. The \emph{star arboricity} of a graph is the size of the minimum forest, each tree of which is a star -- a star is a tree with at most one non-leaf node -- into which the edges of the graph can be partitioned. 
	There are similar notions like \emph{anarboricity}, \emph{pseudoarboricity} and \emph{thickness}~\cite{A1,A2,Harary,Alon-arboricity} for an undirected graph. The common theme running through all such problems is about figuring out the minimum number of structures, e.g., trees, paths, stars, etc., needed to form a disjoint cover of the edge set of the graph.   
}

\medskip
The notion of covering a point set with different geometric structures have been studied in the literature~\cite{brass, keszegh2013covering, dillencourt2004geometric, dumitrescu2014covering, grantson2006covering}. The common theme running through all such problems is about figuring out the minimum number of structures, e.g., trees, paths, line segments, etc., needed to form a cover of the point set.  
Given a set of points, a {\it covering path} is a polygonal path that visits all the points and similarly a {\it covering tree} is a tree whose edges are line segments that jointly cover all the points. There has been considerable interest in the study of covering paths and trees in recent years, see e.g. \cite{arkin2003, dumitrescu2014}. Covering paths and trees for planar grids have been studied in \cite{keszegh2013covering}, where bounds on the minimum number of line segments of such paths and trees are given. Analogous questions on covering paths and trees for higher dimensional grids have been studied in \cite{dumitrescu2014covering}. Moreover covering paths and trees possessing the additional property of being non-crossing have been studied in \cite{biniaz2024, cerny2007}. Given a set $S$ of $n$ points in the plane, the problem of finding the smallest number $l$ of straight lines needed to cover all $n$ points in $S$ have been studied in \cite{grantson2006covering}, where bounds on the time complexity of this problem in terms of $n$ and $l$ (assuming $l$ to be small) is given (see also \cite{chen2021}). In \cite{dillencourt2004geometric}, the notion of {\it geometric thickness} of complete graphs is studied, where geometric thickness of a graph is defined to be the smallest number of layers such that we can draw the graph in the plane with straight-line edges and assign each edge to a layer so that no two edges on the same layer cross. The intersection of a convex body with a lattice is called a {\it convex set of lattice points}. Several problems, conjectures and results on covering a convex set of lattice points by minimum number of lines, hyperplanes, or other subspaces have been discussed in \cite{brass}.

\medskip

On the other hand, incidence theory in geometry~\cite{books/daglib/0018467,books/daglib/0080838} studies questions about finding the maximum possible number of pairs $(p,\ell)$ such that $p$ is a point belonging to a set of points and $\ell$ is a line belonging to a set of lines and $p$ lies on $\ell$. \remove{We will henceforth say that $\ell$ \emph{covers} $p$.} Incidence between points and other geometric structures like circles, planes, algebraic curves, etc.~have also been studied and bounds on the number of incidences were obtained (see e.g.~\cite{books/daglib/0018467,books/daglib/0080838}). On the other hand, researchers have studied the problems of \emph{point line cover}, or its more general form of \emph{point curve cover}~\cite{KratschPR16,BoissonnatDGK18,AfshaniBDN16,LangermanM05}.  These problems consist of a set $P$ of $n$ points on the plane and a positive integer $k$, and the question is whether there exists a set of at most $k$ lines/hyperplanes/curves which cover all points in $P$. They are computationally hard problems, motivated from {\sc set cover}, and the effort has been mostly in parametrized complexity where researchers focused on finding tight kernels~\cite{saketbook15} for the problems~\cite{KratschPR16,BoissonnatDGK18,AfshaniBDN16,LangermanM05}.    

One of the main problems we investigate in this paper is the geometric covering number of a planar grid for orthoconvex curves (a rectilinear curve is called orthoconvex if for any two points
with the same x/y coordinate lying inside the curve, the
vertical/horizontal line segment joining the two points also lies
inside the curve). Despite being natural, to the best of our knowledge, this problem has not been looked at in the literature and our work seems to be the first to address it. To us, the general problem seems difficult and we currently do not have any tight bounds for it. Note that, the number of such curves is exponential in the grid-size, which might be the cause of hardness of this problem. So, as a first step towards the attack on this problem, we restrict our attention to the simplest orthoconvex curves. These are the curves that have few (constant) number of non-convex vertices (which we call inner-corners). Note that, the number of such curves is polynomial in the grid-size. We prove the following bounds (of which, the first one is tight) on the geometric covering number of a square grid for orthoconvex curves with at most $1$ and $2$ inner-corners. Section~\ref{sec:closedcurve} has the following results.
\begin{itemize}
	\item If $m$ orthoconvex curves with at most one inner corner cover the $n\times n$ grid, then $m \geq 2n/5$. 
	\item We need at least $2n/7$ orthoconvex curves with at most two inner corners to cover an $n\times n$ grid. 
\end{itemize}
\remove{
	\newtheorem*{thm:orthoconv}{Theorem \ref{thm:orthoconv}}
	\begin{thm:orthoconv}
		If $m$ orthoconvex curves with at most one inner corner cover the $n\times n$ grid, then $m \geq 2n/5$.
	\end{thm:orthoconv}
	
	\newtheorem*{thm:2corner}{Theorem \ref{thm:2corner}}
	\begin{thm:2corner}
		We need at least $2n/7$ orthoconvex curves with at most two inner corners to cover an $n\times n$ grid
	\end{thm:2corner}
}
Next we mention the following result that is a converse to the line-covering problem, which also gives a counter-example to a conjecture recently published in \cite{sheffer_2022}. This result is discussed in Section~\ref{sec:linecover}.
\begin{itemize}
	\item There exists a finite set $P$ of $n^2$ points in $\mathbb{R}^2$ which can be covered with $n$ lines but no subset of $P$ of size $\Omega (n^2)$ can be contained in a projective transformation of a rectangular grid of size $o(n^3)$.
\end{itemize}
\remove{
	\newtheorem*{thm:line-converse}{Theorem \ref{thm:line-converse}}
	\begin{thm:line-converse}
		There exists a finite set $P$ of $n^2$ points in $\mathbb{R}^2$ which can be covered with $n$ lines but no subset of $P$ of size $\Omega (n^2)$ can be contained in a projective transformation of a rectangular grid of size $o(n^3)$.
	\end{thm:line-converse}
}
In addition to the above major results, we also prove some minor ones on the geometric covering number for lines, algebraic curves, monotonic curves and convex curves. These results are relatively immediate, but our paper offers a full collection and puts them into a coherent picture.

\medskip

\paragraph*{Notations:} $[x]$ will denote the set of natural numbers $\{1, 2, \ldots, x\}$. $P$ will denote a set of $n$ points in dimension $\R^d$. Unless otherwise stated, $P$ will be finite. 

\medskip

\paragraph*{Organization of the paper:} In this paper, we study the notion of geometric covering number for a few types of curves. For most of the cases, our point set is a grid that we want to cover with a particular kind of curve. We start with lines, the simplest curve, covering a finite grid in Section~\ref{sec:linecover}. We also investigate a converse question of covering in Section~\ref{ssec:conversecover}. Very simply put, the converse question deals with the following notion -- if there is a guarantee that some lines cover an ``unknown'' point set, then can we say something about the configuration of the point set? From lines, we move onto monotone curves in Section~\ref{sec:monotone-cover}. Section~\ref{sec:closedcurve} considers three types of closed curves -- circles, convex curves, orthoconvex curves. \remove{Section~\ref{sec:non-congruent} brings into focus covering with non-congruent curves. In Section~\ref{sec:smallcurve}, we study covering by some small curves.}Finally, Section~\ref{sec:conclude} sums up the findings in this work and mentions some open problems. We feel that our work will motivate studying the \emph{geometric covering number} for other point configurations and other types of curves.

\section{Covering by lines and its converse problem}\label{sec:linecover}
\noindent	
In the first part of this section, we consider covering grids by lines\remove{\footnote{We mean lines clipped by the grid.}}\remove{(the bounds can be easily obtained; we have included it for the sake of completeness)}. In the next part, we consider a ``converse'' question -- if a set of $n^2$ points in $\mathbb{R}^2$ is covered by $n$ lines, then can we say something about the configuration of the points? 
\subsection{Covering by lines}\label{ssec:linecover}
\noindent
Note that for any two points there exists a line covering them. Therefore, \ga$ \leq \frac{\size{P}}{2}$ (the equality is achieved for any set of points in general position). Now let $\ell (P)$ denote the maximum number of points in $P$ any line can cover. Then we have \ga$\geq \size{P}/\ell (P)$. Therefore, we get 
\begin{equation}
	\frac{\size{P}}{\ell (P)} \leq \mbox{\ga} \leq \frac{\size{P}}{2}.
	\label{eq:linebound}
\end{equation}
Now we consider the case when $P = [k_1]\times \cdots \times [k_d]$. Let $M:=\max \{ k_1,\ldots ,k_d\}$. We prove the following:
\begin{prop}
	$\ell (P)=M.$
	\label{prop:line}
\end{prop}
\begin{proof}
	First we show that $\ell (P)\leq M$ by induction on $d$. The base case $d=1$ is obvious. Now we proceed to the induction step. Let $L$ be a line segment that lies inside the rectangular parallelepiped $[1,k_1]\times\cdots\times [1,k_d]$. Then $L$ has length at most $\sqrt{\sum_{i=1}^d (k_i-1)^2}.$ Now let $x:=(x_1,\ldots ,x_d)$ and $y:=(y_1,\ldots ,y_d)$ be two distinct points of $P$ lying on $L$. If $x_i=y_i$ for some $i$, then $L$ lies inside a lower dimensional rectangular parallelepiped. Then, by induction hypothesis, $L$ covers at most $\max \{k_j\mid j\ne i\}\leq M$ many points. So let us assume $x_i\ne y_i$ for all $i=1,\ldots ,d$. Then the distance between $x$ and $y$ is at least $\sqrt{d}$. Suppose $L$ covers a total of $t$ points of $P$. Then we have 
	$$ (t-1)\sqrt{d}\; \leq\; \sqrt{\sum_{i=1}^d (k_i-1)^2}\; \leq\; \sqrt{d}\cdot\max\{k_1-1,\ldots ,k_d-1\}  $$
	and this implies $t\leq\max\{k_1,\ldots ,k_d\}$. Therefore, we conclude that $\ell (P)\leq M$. On the other hand, there clearly exist lines covering $M$ points, namely the lines parallel to the coordinate axis $i_0$, where $M=k_{i_0}$. Hence, we have shown that $\ell (P)=M$.
\end{proof}

\medskip	

So, Proposition~\ref{prop:line} implies that 
$$\mathcal{G}_C(P) \geq \frac{\prod_{i=1}^dk_i}{M} \geq \min\left\{ \prod_{i \ne 1} k_i, \ldots,\prod_{i\ne d}k_i \right\}:=N.$$
On the other hand, \ga$\leq N$ since there clearly exists an explicit covering of $P$ by $N$ lines (namely, by the lines parallel to the coordinate axis $i_0$). Therefore, we get that
$$\mathcal{G}_C(P)= \min\left\{ \prod_{i\ne 1}k_i,\ldots ,\prod_{i\ne d}k_i \right\}.$$ Next we mention the following two special cases.
\begin{enumerate}
	\item[(i)] For $P=[k_1]\times [k_2]$ ($2\times 2$ grid),
	$\mathcal{G}_C(P)=\min\{k_1,k_2\}.$
	
	\item[(ii)] For $P=\{0,1\}^d$ (hypercube),
	$\mathcal{G}_C(P)=2^{d-1}.$
	Note that for hypercube both inequalities of Equation~\ref{eq:linebound} become tight.
\end{enumerate}

\begin{rmk}[Skew lines]
	We say that a line is {\it skew} if it is not parallel to x or y-axis. We look at the question of covering an $n\times n$ grid by the minimum number of skew lines.
	
	\medskip		
	
	Note that the boundary of the $n\times n$ grid contains $4n-4$ points. Now any skew line can contain at most 2 points from the boundary. So we need at least $2n-2$ skew lines to cover the grid. Also note that the $n\times n$ grid can be covered by $2n-2$ skew lines (consider the $2n-3$ lines parallel to the off-diagonal except the ones which pass through the bottom-left and top-right corners and these two corners are covered by the main diagonal).
	
	\medskip
	
	It is an {\it open problem} to find the minimum number of skew hyperplanes required to cover the $d$-dimensional hypercube. Currently, the best known lower bound for the above problem is $d/2+1$ (see \cite{ivanisvili2025, saurmann2024, linial2005})
	
\end{rmk}

\subsection{On the converse of the covering problem}\label{ssec:conversecover}
\noindent
We have seen that an $n\times n$ grid can be covered by $n$ lines. Here we look at the converse question, namely, if a set of $n^2$ points in $\mathbb{R}^2$ is covered by $n$ lines then can we say something about the configuration of the points?

\medskip	

Suppose a set of $n^2$ points is covered by $n$ lines. Then there exists a line containing $\Omega (n)$ points, since otherwise the total number of points is less than $n^2$. Now if this line contains $o(n^2)$ points, then there exists another line containing $\Omega (n)$ points. By continuing this, we can say that there exists a set of lines each containing $\Omega (n)$ points such that the total number of points in the union of these lines is $\Theta (n^2)$.

\medskip	

Now the following question seems natural. If a set $P$ of $n^2$ points is covered by $n$ lines, then does there always exist a subset of $P$ of size $\Theta (n^2)$  which can be put inside a grid of size $\Theta (n^2)$, possibly after applying a projective transformation? We show that the answer is no.

\medskip	

\begin{thm}
	\label{thm:line-converse}
	There exists a finite set $P$ of $n^2$ points in $\mathbb{R}^2$ which can be covered with $n$ lines but no subset of $P$ of size $\Omega (n^2)$ can be contained in a projective transformation of a rectangular grid of size $o(n^3)$.
\end{thm}

\begin{proof}		
	Given any two distinct points $p,p'\in\mathbb{R}^2$, we denote by $\ell(p,p')$ the unique line in $\mathbb{R}^2$ that contains both $p$ and $p'$.
	By an $s\times t$ grid, we mean a point set that can be obtained by a projective transform $f$ of the set $[t]\times [s]$. By a ``horizontal line'' of the grid, we mean a line $\ell(f(1,j),f(t,j))$ for some $j\in [s]$, and by a ``vertical line'' of the grid, we mean a line $\ell(f(i,1),f(i,s))$, for some $i\in [t]$. The ``size'' of an $s\times t$ grid is $st$, i.e., the number of points in it. Note that every horizontal line of a grid intersects every vertical line of the grid (since there is a point of the grid that is contained in both of them).
	
	\medskip	
	
	For each $i\in [n]$, let $L_i$ denote the line with equation $y=i$ and let $\mathcal{L}=\{L_i\}_{1\leq i\leq n}$.
	Let $\mathcal{P}$ be the set of points defined as follows. Define $P_1$ to be some set of $n$ distinct points from the line $L_1$. For each $1<i\leq n$, we define $P_i$ to be a set of $n$ distinct points from $L_i$ that do not lie on any of the lines formed by points on other lines, i.e. in $\{\ell(p,p')\colon p\neq p'$ and $p,p'\in\bigcup_{1\leq j\leq i-1} P_j\}$. Let $\mathcal{P}=\bigcup_{1\leq i\leq n} P_i$. Let $m=|\mathcal{P}|$. Note that we have $|\mathcal{L}|=n$. We claim that for any $\mathcal{P}'\subseteq\mathcal{P}$ such that $|\mathcal{P'}|=\Omega(m)=\Omega(n^2)$, any grid that contains all the points of $\mathcal{P}'$ has size $\Omega(n^3)$. 
	
	\medskip		
	
	Note that by our construction, if any line contains two points $p,p'\in\mathcal{P}$ such that $p\in P_i$ and $p'\in P_j$, where $i\neq j$, then $p$ and $p'$ are the only points in $\mathcal{P}$ that are contained in that line. This implies that the following property is satisfied by $\mathcal{P}$ and $\mathcal{L}$.
	\medskip
	
	\noindent (*) Any line in $\mathbb{R}^2$ that contains more than two points in $\mathcal{P}$ belongs to $\mathcal{L}$.
	\medskip
	
	Since every line in $\mathcal{L}$ contains exactly $n$ points of $\mathcal{P}$, we then have another property.
	\medskip
	
	\noindent (+) Any line in $\mathbb{R}^2$ contains at most $n$ points in $\mathcal{P}$.
	\medskip
	
	Let $\mathcal{P}'\subseteq\mathcal{P}$ be such that $|\mathcal{P}'|=\Omega(n^2)$. Consider any grid $\mathbb{G}$ that contains all the points of $\mathcal{P}'$. Let $\mathbb{G}$ be an $s\times t$ grid. Let $h_1,h_2,\ldots,h_s$ denote the horizontal lines of $\mathbb{G}$ and let $v_1,v_2,\ldots,v_t$ denote the vertical lines of $\mathbb{G}$. Suppose for the sake of contradiction that there exist $i\in [s]$ and $j\in [t]$ such that both the lines $h_i$ and $v_j$ contain at least 3 points of $\mathcal{P}$ each. Then by property~(*), $h_i$ and $v_j$ are both lines in $\mathcal{L}$. But as $h_i$ and $v_j$ intersect, they are two lines in $\mathcal{L}$ that intersect, which is a contradiction, since the lines in $\mathcal{L}$ are all parallel to each other (Note that parallel lines under a projective transformation may not be parallel but they do not intersect at any of the $s\times t$ grid points defined.). Thus, we can conclude without loss of generality that for each $i\in [s]$, the horizontal line $h_i$ of $\mathbb{G}$ contains at most two points from $\mathcal{P}$, and hence at most two points from $\mathcal{P}'$. Since every point in $\mathcal{P}'$ is contained in at least one horizontal line of $\mathbb{G}$, we have that $s\geq |\mathcal{P}'|/2$ and therefore $s=\Omega(n^2)$. By property~(+), each vertical line of $\mathbb{G}$ can contain at most $n$ points of $\mathcal{P}'$, and therefore, $t\geq |\mathcal{P}'|/n$, which implies that $t=\Omega(n)$. Thus the size of the grid $\mathbb{G}$ is $st=\Omega(n^3)$.
\end{proof}

\subsubsection*{Regarding a conjecture on incidence geometry}

Given a set of points $P$ and a set of lines $L$ in the plane, an {\it incidence} is a pair $(p,l)\in P\times L$ such that $p\in l$. Let $I(P,L)$denote the number of incidences. By a theorem of Szemer\'{e}di and Trotter~\cite{szemeredi1983}, it  is known that $I(P,L)=O(m^{2/3}n^{2/3}+m+n)$. On the other hand, Erd\"{o}s~\cite{erdos1985} gave construction of $m$ points and $n$ lines with $\Theta (m^{2/3}n^{2/3}+m+n)$ incidences. In this construction, the point set is a rectangular section of the integer lattice. One can also obtain somewhat different point sets by applying various projective transformations. In \cite{sheffer_2022}, the following conjecture was mentioned.

\begin{conj}
	Consider sufficiently large positive integers $m$ and $n$ that satisfy $m = O(n^2)$ and $m = \Omega(\sqrt{n})$. Let $P$ be a set of $m$ points and $L$ be a set of $n$ lines, both in $\R^2$, such that $I(P,L) = \Theta (m^{2/3} n^{2/3})$. Then there exists a subset $P' \subset P$ such that $\size{P'} =  \Theta(m)$ and $P'$ is contained in a section of the integer lattice of size $\Theta(m)$, possibly after applying a projective transformation to it.	
\end{conj}

Note that, our construction (with with $n^2$ points, $n$ lines and $n^2$ incidences) in the proof of Theorem~\ref{thm:line-converse} provides a counter-example to the above conjecture. This was communicated to Prof. Adam Sheffer who told us that this only exposes a typo in the statement of the conjecture which is more interesting and challenging when $m=o(n^2)$. Note that in our construction we have $m=n^2$. Below we state the modified conjecture.

\begin{conj}
	Consider sufficiently large positive integers $m$ and $n$ that satisfy $m = o(n^2)$ and $m = \Omega(\sqrt{n})$. Let $P$ be a set of $m$ points and $L$ be a set of $n$ lines, both in $\R^2$, such that $I(P,L) = \Theta (m^{2/3} n^{2/3})$. Then there exists a subset $P' \subset P$ such that $\size{P'} =  \Theta(m)$ and $P'$ is contained in a section of the integer lattice of size $\Theta(m)$, possibly after applying a projective transformation to it.	
\end{conj}

\subsection{Covering by algebraic curves}

In this subsection, we address the question of covering a grid by algebraic curves. The answer comes as a direct application of the famous Combinatorial Nullstellensatz Theorem due to Noga Alon.

\begin{lem}[Combinatorial Nullstellensatz~\cite{alon_1999}]
	Let $f = f(x_1,\ldots , x_d)$ be a polynomial in $\mathbb{R}[x_1,\ldots ,x_d]$.
	Suppose the degree $\deg(f)$ of $f$ is $\sum_{i=1}^d t_i$ where each $t_i$ is a non-negative integer, and suppose the coefficient of $\prod_{i=1}^d {x_i}^{t_i}$ in $f$ is non-zero. Then, if $S_1,\ldots ,S_n$ are subsets of $\mathbb{R}$ with $|S_i| > t_i$, there are $s_1\in S_1, s_2 \in S_2,\ldots , s_d \in S_d$ so that $f(s_1,\ldots ,s_d)\ne 0$.
	
	\label{lem:null}
\end{lem}

\begin{thm}
	Suppose the $n\times n$ grid is covered by $m$ algebraic curves of degree at most $k$. Then $m\geq n/k$.
\end{thm}

\begin{proof}
	Suppose $m<n/k$. Let the algebraic curves defined by $f_1(x,y)=0,\ldots ,f_m(x,y)=0$ cover the $n\times n$ grid, where $\deg(f_i)\leq k$. Then the polynomial $f(x,y):=\prod_{i=1}^m f_i(x,y)$ vanishes at each grid point. Suppose $\deg (f) = t_1+t_2$ with the coefficient of $x^{t_1}y^{t_2}$ in $f$ being non-zero. Now note that $t_i\leq t_1+t_2 = \deg (f)\leq mk < n$, for each $i=1,2$. So by Lemma~\ref{lem:null}, there exists a grid point $(s_1, s_2)$ so that $f(s_1, s_2) \ne 0$ and we arrive at a contradiction. Therefore, we conclude that $m\geq n/k$.
\end{proof}

\begin{cor}
	\ga $=\lceil n/k\rceil$, where $P$ is an $n\times n$ grid and $C$ denotes algebraic curves of degree at most $k$.
\end{cor}

\begin{proof}
	The lower bound follows from the previous theorem and the upper bound follows from covering by lines and then considering a set of $k$ lines as one curve of degree $k$.
\end{proof}

\begin{rmk}[Irreducible algebraic curves]
	By a result of Bombieri and Pila~\cite{bombieri}, an {\it irreducible} algebraic curve of degree $k$ can contain at most $O(n^{1/k})$ points from an $n \times n$ grid and hence, the minimum number of irreducible algebraic curves of degree $k$ to cover the $n\times n$ grid is at least $\Omega (n^{2-1/k})$.
	
\end{rmk}

Using the same reasoning as in the previous theorem and corollary, one also has the following result on covering the $n_1\times\cdots \times n_d$ grid by algebraic hypersurfaces.

\begin{thm}
	The minimum number of algebraic hypersurfaces of degree at most $k$ needed to cover the $n_1\times\cdots \times n_d$ grid is equal to $\lceil n/k\rceil$, i.e., \ga $=\lceil n/k\rceil$, where $P$ is an $n_1\times\cdots \times n_d$ grid and $C$ denotes algebraic hypersurfaces of degree at most $k$.
\end{thm}

\section{Covering by monotonic curves}\label{sec:monotone-cover}
\noindent
In this section, we consider the case when the the curve is \emph{monotonic}.
\begin{defn}
	Let $f:[0,1]\to\mathbb{R}^d$ be a curve and suppose $f(t)=(f_1(t),\ldots ,f_d(t))$ for $t\in [0,1]$. Then $f$ is called \textbf{monotonic} if it satisfies the following property: 
	$t_1\leq t_2 \; \Rightarrow \; f_i(t_1)\leq f_i(t_2) \text{ for    each } i=1,\ldots ,d.$ 
\end{defn}

\paragraph*{Preliminaries on posets:} Let $P$ be a finite set and $\leq$ be a binary relation on $P$. Then $(P,\leq )$ is called a {\it partially ordered set} ({\it poset} for short) if it satisfies the following conditions.
\begin{enumerate}
	\item Reflexivity: $x\leq x$ for all $x\in P$.
	
	\item Antisymmetry: For $x,y\in P$, $x\leq y$ and $y\leq x$ implies $x=y$.
	
	\item Transitivity: For $x,y,z\in P$, $x\leq y$ and $y\leq z$ implies $x\leq z$.
\end{enumerate}

A subset $S$ of $P$ is called a {\it chain} if for every pair $x,y\in S$, either $x\leq y$ or $y\leq x$ (i.e. any two elements of $S$ are comparable). Similarly a subset $T$ of $P$ is called an {\it antichain} if for every pair $x,y\in T$, neither $x\leq y$ nor $y\leq x$ (i.e. no two elements of $T$ are comparable).

We say that $x<y$ if $x\leq y$ and $x\ne y$. A poset $(P,\leq )$ is called {\it graded} if $P$ is equipped with a rank function $r:P\to\mathbb{N}$ which satisfies the following two properties: (i) If $x < y$ then $r(x)<r(y)$, and (ii) If $x<y$ and there does not exist any $z$ such that $x<z<y$, then $r(y)=r(x)+1$.

Let $(P,\leq)$ be a graded poset with rank function $r$. Let $P_i$ denote the set of elements of $P$ of rank $i$, i.e. $P_i=\{x\in P: r(x)=i\}$. Then $P=\biguplus_{i=1}^n P_i$, where $n=\max_{x\in P}{r(x)}$. The poset $P$ is called {\it rank symmetric} if $|P_i|=|P_{n-i}|$ for all $i$. It is called {\it rank unimodal} if there is some $k$ such that $|P_1|\leq\cdots\leq |P_k|\geq\cdots\geq |P_n|$. Let $m:=\max_{i}|P_i|$ be the size of the largest rank level. $P$ is called {\it Sperner} if no antichain within it is larger than $m$. 

\medskip

Given a finite subset $P$ of $\mathbb{R}^d$, we define the poset $\mathcal{P} := (P,\leq)$ as follows. For $x:=(x_1,\ldots , x_d)\in\mathbb{R}^d$ and $y:=(y_1,\ldots ,y_d)\in\mathbb{R}^d$, we define $x\leq y$ if $x_i\leq y_i$ for  $i=1,\ldots ,d$. We prove the following proposition.
\begin{prop}
	Let $w(\mathcal{P})$ denote the size of the largest antichain, called the width, of $\mathcal{P}$. Then \ga = $w(\mathcal{P})$, where $P$ is any point set and $C$ denotes monotonic curves.
\end{prop}
\begin{proof}
	Let $x_i\in\mathcal{P}$ for $i=1,\ldots,r$. Then note that $x_1\leq\cdots \leq x_r$ is a chain if and only if $x_1,\ldots ,x_r$ lie on the same curve (which is monotonic). Therefore, \ga~equals the number of chains in the minimal chain decomposition of $\mathcal{P}$, which by Dilworth's theorem~\cite{dilworth} equals the size of the largest antichain of $\mathcal{P}$. Hence \ga$=w(\mathcal{P})$.
\end{proof}

Note that the poset $\mathcal{P}$ can be decomposed into $w(\mathcal{P})$ many disjoint chains. Therefore, the points in $P$ can be covered by \ga~many monotonic curves such that  \textit{no two curves intersect at a point of $P$}.

\medskip

Now we consider the special case when $P=[k_1] \times \cdots \times [k_d]$, where $k_i \in \mathbb{N}$ for each $i$. It is known that the poset $\mathcal{P}$ is graded and Sperner. In fact, $\mathcal{P}$ is a {\it Peck poset} (i.e., rank symmetric, rank unimodal and Sperner) since it is a product of Peck posets, namely chains (see Theorem 6.2.1 of \cite{engel} for a proof of the fact that product of Peck posets is Peck). Therefore, $w(\mathcal{P})$ equals the size of the largest rank of $\mathcal{P}$, i.e., 
$w(\mathcal{P})=\max_m A_m = A_{\lfloor (k_1+\cdots +k_d+d)/2\rfloor},$
where, $A_m$ equals the number of solutions of the equation $x_1+\cdots +x_d=m$ such that $x_i\in[k_i]$ for each $i=1,\ldots, d$. We mention the following two special cases.
\begin{enumerate}
	
	\item[(i)] When $d=2$ (i.e., $P$ is a $2\times 2$ grid), $w(\mathcal{P})=\max_m A_m = \min\{k_1,k_2\}$.
	
	\item[(ii)] When $k_1= \cdots= k_d = 2$ (i.e., $\mathcal{P}$ is a Boolean lattice), $w(\mathcal{P})=\max_m A_m = A_{\left\lfloor d + \frac{d}{2}\right\rfloor}=\binom{d}{\lfloor d/2\rfloor}$.		
\end{enumerate}

\section{Covering by closed curves}\label{sec:closedcurve}
\noindent
In this section, we move onto closed curves; we consider covering grids by circles, convex curves and orthoconvex curves. Notice that the curves need not be of the same size, e.g., when we are considering covering by circles, all the circles need not be of the same size.
\subsection{Covering by circles} 
A circle contains at most $O(n^{\epsilon})$ points from an $n\times n$ grid for every $\epsilon >0$ (see e.g. \cite {guth}). Therefore, the minimum number of circles required to cover an $n \times n$ grid is $\Omega (n^{2-\epsilon})$, for every $\epsilon >0$.

\medskip	

Regarding upper bound, note that there is a covering of the $n\times n$ grid by $O(n^2/\sqrt{\log n})$ circles. This is obtained by choosing a corner of the grid and drawing all concentric circles such that each of them is incident to at least one grid point. The number of such circles is $O(n^2/\sqrt{\log n})$ by a well known theorem of Ramanujan and Landau (\cite{ramanujan},\cite{landau}). The theorem says that the number of positive integers that are smaller than $n$ that are the sum of two squares is $\Theta (n/\sqrt{\log n})$.

\medskip

We sum it up as the following.

\begin{prop}
	$\Omega (n^{2-\epsilon}) \leq $ \ga $ \leq O(n^2/\sqrt{\log n})$, where $P$ is an $n\times n$ grid and $C$ denotes circles.
\end{prop}

\subsection{Covering by convex curves}

A closed convex curve intersects non-trivially with a horizontal grid line if it contains more than two points from the line. Note that, any closed convex curve can intersect at most two horizontal grid lines non-trivially. This follows from the following lemma.
\begin{lem}
	If a closed convex curve intersects a horizontal grid line non-trivially, then it must lie entirely on one side of that line.
	\label{lem:oneside}
\end{lem}
\begin{proof}
	Suppose the curve intersects a horizontal line at three points $p, q, r$, where $q$ lies in the interior of line segment $[p,r]$. Since the curve is convex, there exists a line $L$ through $q$ such that the curve lies entirely on one side of $L$ (hyperplane separation theorem). Now if $L$ is different from the horizontal line, then $p$ and $r$ lie on different sides of $L$. But since the curve lies on one side of $L$, it can not pass through both $p$ and $r$, a contradiction. Therefore, $L$ is same as the horizontal line and the curve lies entirely on one side of this line.
\end{proof}
\begin{thm}
	The points of the $n\times n$ grid cannot be covered with less than $n/2$ closed convex curves, i.e. \ga $ \geq n/2$ where $P$ is an $n\times n$ grid and $C$ denotes closed convex curves.
\end{thm}
\begin{proof}
	Suppose, for the sake of contradiction, that $C_1,C_2,\ldots ,C_k$ are $k$ closed convex curves such that they together cover every point of the $n\times n$ grid and that $k < n/2$. Now, there are $n$ horizontal grid lines. By Lemma~\ref{lem:oneside} above, each $C_i$ can have a non-trivial intersection with at most 2 horizontal grid lines. So we conclude that there is some horizontal grid line such that no curve in $C_1, C_2,\ldots ,C_k$ has a non-trivial intersection with that line. Now consider the points on that horizontal line. There are $n$ points on this line. Each curve in $C_1, C_2, \ldots ,C_k$ can cover at most two points from that line, since none of them intersects non-trivially with this horizontal line. But then, since $k < n/2$, there must be some point on this horizontal line that is not covered by any curve in  $C_1, C_2,\ldots, C_k$, which is a contradiction.
\end{proof}

Almost same argument can be used to show that at least $\min{\{\lceil m/2\rceil , \lceil n/2\rceil \}}$ convex curves are required to cover an $m \times n$ grid.

\begin{defn}
	In $\mathbb{R}^d$, we say that a closed convex hypersurface intersects a hyperplane non-trivially if it intersects the hyperplane in at least $d+1$ points such that one of these $d+1$ points lie in the interior of the convex hull of the rest of the $n$ points.
\end{defn}
With this definition, the same argument (as in the 2-dimensional case) will go through. So, finally we have the following theorem by inductive argument (where induction is on the dimension $d$ of the grid).
\begin{thm}
	The minimum number of closed convex hypersurfaces required to cover the $k_1\times\cdots\times k_d$
	grid is $\min{\{\lceil k_1/2\rceil \ldots ,\lceil k_d/2\rceil \}}$, i.e. \ga $ = \min{\{\lceil k_1/2\rceil \ldots ,\lceil k_d/2\rceil \}}$ where $P$ is a $k_1\times\cdots\times k_d$ grid and $C$ denotes closed convex hypersurfaces.
\end{thm}

\begin{rmk}[Strictly convex curves]
	Any strictly convex curve (convex curve which does not contain a line segment) can contain $O(n^{2/3})$ points of an $n\times n$ grid by a theorem of Andrews (\cite{andrews}). Therefore, we need $\Omega (n^{4/3})$ strictly convex curves to cover an $n\times n$ grid. On the other hand, it follows from the result proved in \cite{har} that one can cover the $n\times n$ grid by $O(n^{4/3})$ strictly convex curves. Hence we conclude that \ga $ =\Theta (n^{4/3})$, where $P$ is an $n\times n$ grid and $C$ denotes strictly convex curves.
\end{rmk}

\subsection{Covering by orthoconvex curves}
\begin{figure}[h]
	\centering
	\includegraphics[scale=0.5]{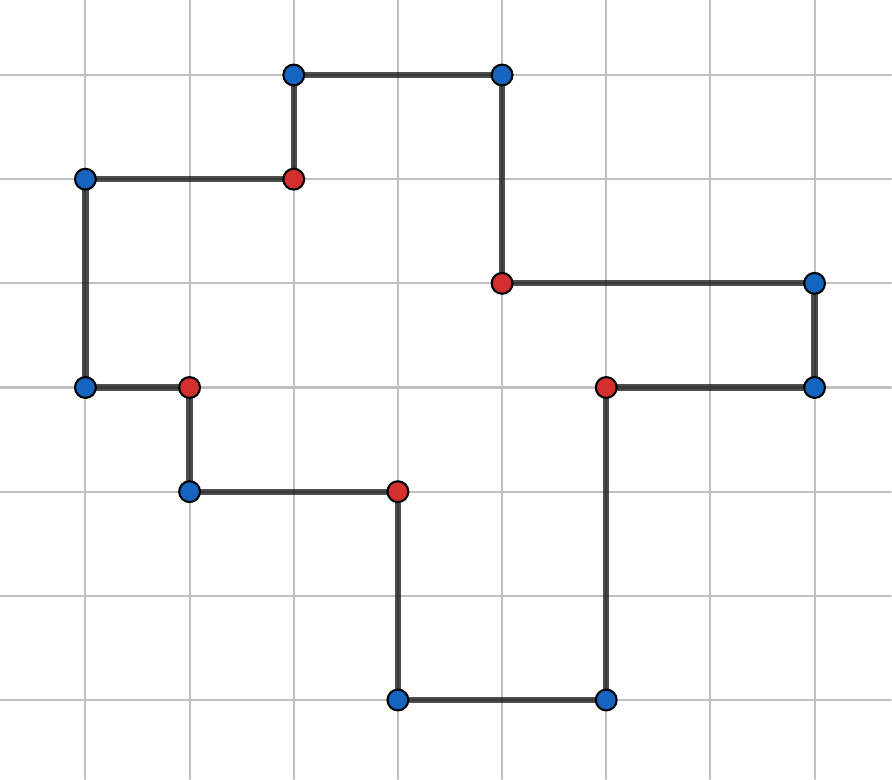}
	\caption{An othoconvex curve and its inner corners (in red)}
	\label{fig:oconv}
\end{figure}

\noindent
A set $K \subseteq \mathbb{R}^2$ is defined to be {\em orthogonally convex} if, for every line $\ell$ that is parallel to one of standard basis vectors $(1,0)$ or $(0,1)$, the intersection of $K$ with $\ell$ is empty, a point, or a single segment. The {\em orthogonal convex hull} of a point set $P\subseteq\mathbb{R}^2$ is the intersection of all connected orthogonally convex supersets of $P$.	If the boundary of orthogonal convex hull (of a set of points) is a simple closed curve then we call it an {\em orthoconvex} curve. An orthoconvex curve has only two types of angles -- 90 and 270 degrees. By {\em inner corner} of an orthoconvex curve, we mean a point where the curve turns by 270 degrees. See Figure~\ref{fig:oconv} for an example of an orthoconvex where the red points are its inner corners.

\medskip	

If an orthoconvex curve (with $k$ inner corners) covers a set of points, then there is also an orthoconvex curve (with $k$ inner corners) covering the same points which is not self-intersecting and all the corners are grid points. This can be done by pushing the sides/edges of the curve ``outwards'' (instead of inwards which corresponds to taking orthoconvex hull) till we hit a grid line. So w.l.o.g., we may impose the following assumptions of `non-self-intersecting' and `corners are grid points'.

\medskip	

\begin{figure}[h]
	\centering
	\includegraphics[scale=0.5]{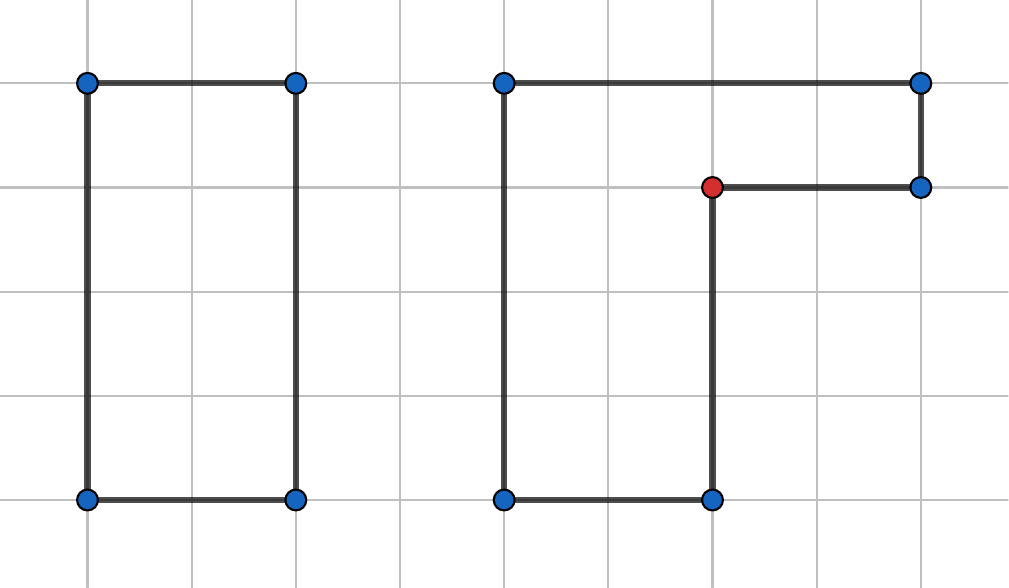}
	\caption{Orthoconvex curves with at most one inner corner}
	\label{fig:1_in_cor}
\end{figure}

In the following, by {\em curve}, we mean an orthoconvex curve having at most one inner corner (Figure~\ref{fig:1_in_cor} shows examples of such curves). We say that a curve {\em hits} a (horizontal or vertical) grid line if the curve has a non-trivial intersection with that grid line (i.e., the curve follows that grid line for some distance, rather than just crossing it).
We say that a collection of curves $C$ {\em hits} a (horizontal or vertical) grid line if there is some curve in $C$ that hits that grid line.
Given a collection of curves $C$, we say that a grid point is {\em exposed} (by $C$) if the grid point is not covered by any curve in $C$, but it lies on a horizontal grid line and a vertical grid line both of which are hit by $C$.
Given a collection of curves $C$, a {\em corner} of $C$ is a corner of the (minimum size) bounding box of $C$. So every collection $C$ of curves has exactly 4 corners.
If a corner of $C$ is an exposed grid point, then we call it an {\em exposed corner}.
We say that a sequence of curves $c_1, c_2,\ldots , c_t$ is {\em good} if for every $i \in \{ 2, 3,\ldots , t \}$, $c_i$ hits a grid line that is hit by $\{ c_1, c_2,\ldots , c_{i-1} \}$.
Clearly, every prefix of a good sequence is also a good sequence.

\begin{lem}
	Let $c_1, c_2,\ldots , c_t$ be a good sequence of curves. Then $\{ c_1, c_2, \ldots , c_t \}$ either: (a) hits at most 5t grid lines, or (b) hits $5t+1$ grid lines and has an exposed corner.
	\label{lem:orthoconvex}
\end{lem}
\begin{proof}
	We prove this by induction on $t$. It is not difficult to see that the lemma is true when $t = 1$. Let $i > 1$ and suppose that the lemma is true for the good sequence $c_1, c_2,\ldots , c_{i-1}$. Let $C = \{ c_1, c_2, \ldots , c_{i-1} \}$. Then either $C$ hits (a) at most $5i-5$ grid lines, or (b) hits $5i-4$ grid lines and has an exposed corner. 
	
	\medskip		
	
	In case (a), since the curve $c_i$ can hit at most 5 grid lines that are not hit by $C$ (recall that $c_i$ hits at least one grid line that is also hit by $C$), we have that $C \cup \{ c_i \}$ can hit at most $5i$ grid lines, and we are done. Next, let us consider case (b). Note that if $c_i$ is a rectangle, then it can hit at most 3 grid lines that are not hit by $C$, and therefore, $C \cup \{ c_i \}$ hits at most $5i - 1$ grid lines, and we are done. So we can assume that $c_i$ is not a rectangle. Also, if there are two grid lines that are hit by both $C$ and $c_i$, then $C \cup \{ c_i \}$ hits at most $5i$ grid lines, and we are done. So we can assume that $c_i$ hits exactly one grid line that is hit by $C$, and therefore, $C \cup \{ c_i \}$ hits exactly $5i+1$ grid lines. In this case, we have to show that one of the corners of $C \cup \{ c_i \}$ is exposed. Let $B$ be the bounding box of $C \cup \{ c_i \}$. Let $g_0, g_1, g_2, g_3$ be the grid lines on which the top, right, bottom, and left borders of $B$ lie. Clearly, each of $g_0, g_1, g_2, g_3$ is hit by either $C$ or $c_i$ or both. Since $c_i$ hits exactly one grid line that is hit by $C$, we have that at most one of $g_0, g_1, g_2, g_3$ is hit by both $C$ and $c_i$. This implies that $C$ and $\{ c_i \}$ do not have shared corners. Note that a corner $v$ of $C$ is exposed, and a corner $v^\prime$ of $\{ c_i \}$ is exposed. If each of $g_0, g_1, g_2, g_3$ is hit by $C$, then $v$ is an exposed corner of $C \cup \{ c_i \}$ (observe that $v$ cannot be covered by $c_i$, because if it is, it has to be a corner of $\{ c_i \}$, which would mean that $C$ and $\{ c_i \}$ have a shared corner) and we are done. Similarly, if each of $g_0, g_1, g_2, g_3$ is hit by $c_i$, then $v^\prime$ is an exposed corner of $C \cup \{ c_i \}$ and we are again done. Thus we can assume that neither $C$ nor $c_i$ hits all the grid lines $g_0, g_1, g_2, g_3$. Recall that all grid lines except at most one in $g_0, g_1, g_2, g_3$ are hit by exactly one of $C$ or $c_i$. Then there exists some $j \in \{0,1,2,3\}$ such that one of $g_j, g_{j+1 \mod 4}$ is hit by $C$ and not by $c_i$, and the other is hit by $c_i$ and not by $C$. Then the grid point that is contained in both the grid lines $g_j$ and $g_{j+1 \mod 4}$ is an exposed corner of $C \cup \{ c_i \}$. This completes the proof.
\end{proof}

\begin{figure}[h]
	\centering
	\includegraphics[scale=0.5]{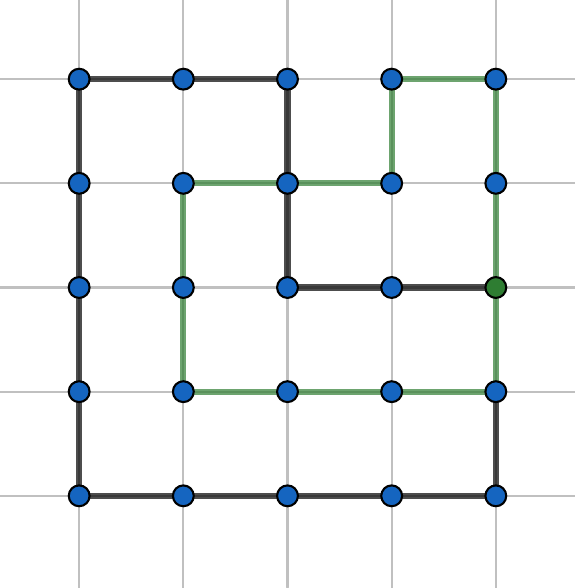}
	\caption{Covering of $5\times 5$ grid by two orthoconvex curves (with at most one inner corner)}
	\label{fig:5by5}
\end{figure}

\begin{thm}
	If $m$ orthoconvex curves with at most one inner corner cover the $n\times n$ grid, then $m \geq 2n/5$.
	\label{thm:orthoconv}
\end{thm}
\begin{proof}
	Let $C$ be a collection of $m$ curves that cover the $n \times n$ grid. For two curves $c$ and $d \in C$, we say that $c R d$ if there is a grid line that is hit by both $c$ and $d$. Let $R^*$ be the transitive closure of $R$. Clearly, $R^*$ is an equivalence relation. Let $S_1, S_2, \ldots , S_p$ be the equivalence classes of $R^*$. We need the following claims for the proof.
	
	\begin{clm} For each $i \in [p]$, $S_i$ does not expose any grid point.
		\label{clm:1}
	\end{clm}
	\begin{proof}
		Suppose for some $i \in [p]$, $S_i$ exposes a grid point $v$. That is, $v$ is not covered by $S_i$, but both the horizontal grid line as well as the vertical grid line that contains $v$ are hit by $S_i$. Since $C$ covers the whole grid, there is a curve $c \in C$ that covers $v$. As $S_i$ does not cover $v$, we have that $c \in C - S_i$. As $c$ covers $v$, $c$ hits either the horizontal grid line containing $v$ or the vertical grid line containing $v$. Since both these grid lines are hit by $S_i$, it follows that there exists some $d \in S_i$ such that $c$ and $d$ hit a common grid line. Then $d R c$, which implies that $c \in S_i$, which is a contradiction. This proves the claim.
	\end{proof}
	\begin{clm}
		The curves of each equivalence class $S_i$ can be arranged in a good sequence.
		\label{clm:2}
	\end{clm}
	\begin{proof}
		Let $G$ be the graph with vertex set $S_i$ and edge set $R$ restricted to $S_i$. We enumerate the curves of $S_i$ in the order in which they are visited by a graph traversal algorithm starting from an arbitrary vertex. Then we get a sequence of the curves in $S_i$ such that before a curve $c$ is encountered in the sequence, we encounter some curve $d$ such that $d R c$ (except for the first curve in the sequence). This sequence is clearly a good sequence of the curves in $S_i$. This proves the claim.
	\end{proof}
	By Lemma~\ref{lem:orthoconvex} and Claims~\ref{clm:1} and \ref{clm:2}, we know that for each $i \in [p]$, $S_i$ hits at most $5 |S_i|$ grid lines. Thus the total number of grid lines that are hit by $C$ is at most $5 ( |S_1| + |S_2| + \cdots + |S_p| ) = 5 |C| = 5m$. If the the curves in $C$ hit $2n$ grid lines, we then have $5m \geq 2n$, which gives $m \geq 2n/5$. Otherwise, suppose that the collection $C$ of $m$ curves, where $m \leq 2n/5$, hits less than $2n$ grid lines. That is, there is some (horizontal or vertical) grid line that is not hit by any curve in $C$. Then every curve in C can cover at most two points on this grid line (if it covers more than two, then the curve hits this grid line). So at most $2m \leq 4n/5$ points on this grid line can be covered by the collection of curves $C$. Hence some points on this grid line are not covered by any curve in $C$, which is a contradiction. So we conclude that $m\geq 2n/5$ and this proves the theorem.
\end{proof}
Note that, the inequality of the above theorem is tight for $n=5$ since the $5\times 5$ grid can be covered by 2 curves (shown in Figure~\ref{fig:5by5}).
As a consequence of the above theorem, we also get the following theorem on orthoconvex curves with {\it at most 2 inner corners}.

\begin{figure}[h]
	\centering
	\includegraphics[scale=0.5]{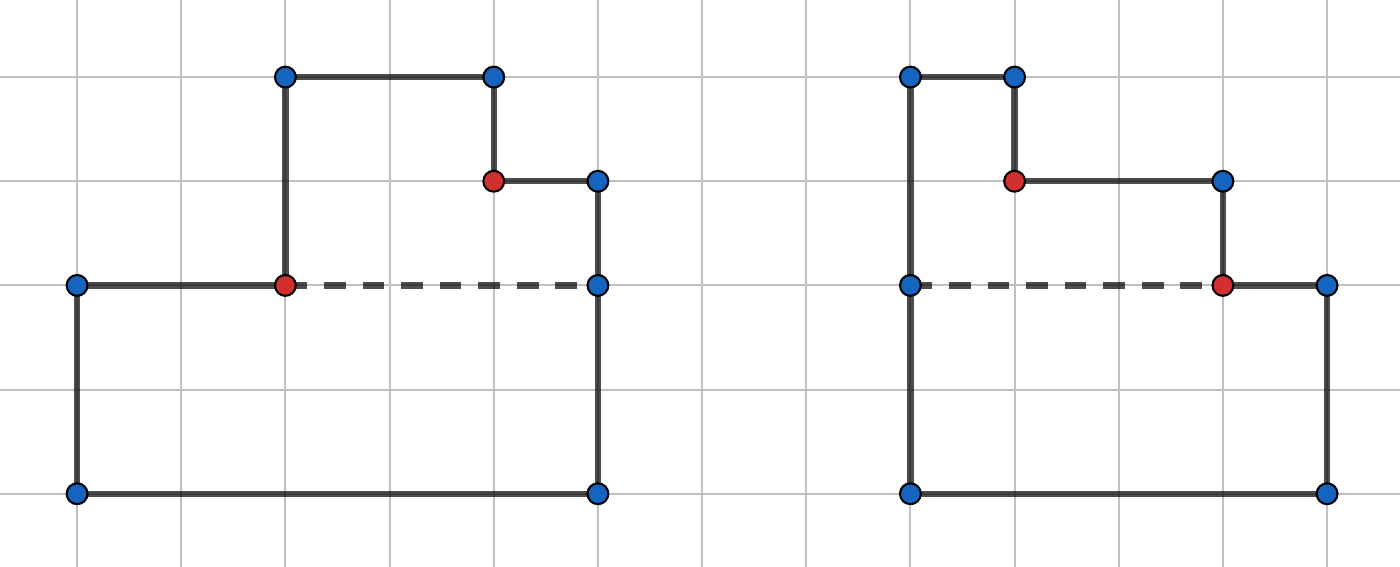}
	\caption{Decomposition of orthoconvex curves with 2 inner corners}
	\label{fig:2in_cor}
\end{figure}

\begin{thm}
	\label{thm:2corner}
	We need at least $2n/7$ orthoconvex curves with at most two inner corners to cover an $n\times n$ grid
\end{thm}
\begin{proof}
	Suppose we have a covering by $m$ such curves. Note that we can decompose each orthoconvex curve with two inner corners into an orthoconvex curve with at most one inner corner and a rectangle (see Figure~\ref{fig:2in_cor}). Hence we obtain a covering by $m$ orthoconvex curves with at most one inner corner and $m$ rectangles. These $m$ orthoconvex curves with at most one inner corner can together hit at most $5m$ grid lines (see proof of Theorem~\ref{thm:orthoconv}). The rectangles together hit at most $2m$ extra grid lines (since each rectangle hit at most two extra grid lines). So the total number of grid lines hit by our original curves is at most $7m$. Since the curves have to hit $2n$ grid lines (by the same reasoning as in proof of Theorem~\ref{thm:orthoconv}), we then have $7m\geq 2n$. Hence, we conclude that $m\geq 2n/7$.	
\end{proof}

\begin{rmk}
	We think that the bound $2n/7$ of Theorem~\ref{thm:2corner} is probably not tight. So a natural problem is to obtain a tight bound for covering by orthoconvex curves with at most 2 inner corners. The next natural follow up question would be: what happens for orthoconvex curves with at most $k$ inner corners for $k=3,4$ etc. It seems our arguments for $k=1,2$ can not be extended to these cases to obtain non-trivial bounds and hence require new ideas. Another question of interest is to find the minimum number of general orthoconvex curves (with no restrictions on the number of inner corners) required to cover an $n\times n$ grid. One can check that for $n=4, \, 5, \, 6, \, 7, \, 8, \, 9 \mbox{ and } 10$, the $n\times n$ grid can be covered by $2, \, 2, \, 2, \, 3, \, 3, \, 3 \mbox{ and } 4$ orthoconvex curves,  respectively. To us, the general problem of orthoconvex curves seems difficult. Note that we have obvious lower and upper bounds of $\lceil (n+1)/4\rceil$ and $\lfloor n/2\rfloor$ respectively, since, any orthoconvex curve can contain at most $4n-4$ grid points (the number of grid pounts on the boundary of an $n \times n$ grid) and on the other hand, an $n \times n$ grid can be covered by $\lfloor n/2\rfloor$ orthoconvex curves. Any improvement over these bounds would be interesting.
\end{rmk}

\section{Conclusion and discussion}
\label{sec:conclude}
In this paper, we mainly discussed the problem of covering a grid (mostly planar) by minimum number of curves of various types. For lines and skew lines, we got the answer by straight forward arguments. An interesting open problem in this direction is to cover the hypercube by minimum number of skew hyperplanes. 
\begin{prob}
	Find the minimum number of skew hyperplanes required to cover the $d$-dimensional hypercube $\{0,1\}^d$.
\end{prob}
For algebraic curves, the answer came as a  consequence of the Combinatorial Nullstellensatz~\cite{alon_1999} but the problem becomes non-trivial if we force the algebraic curves to be irreducible. As a converse to the covering by lines problem, we also show that for a set $P$ of $n^2$ points covered by $n$ lines, it's not true that there always exists a subset of $P$ of size $\Theta(n^2)$ that can be put inside a grid of size $\Theta(n^2)$, possibly after a projective transformation. An open problem here would be to obtain more information about the point configuration $P$. 
\begin{prob}
	If a set of $n^2$ points in the plane is covered by $n$ lines, then what can be said about the point configuration (in additon to Theorem~\ref{thm:line-converse})?
\end{prob}
Next we considered monotonic curves, for which we obtained the answer by applying Dilworth's theorem~\cite{dilworth} on partially ordered sets. Then we looked at three types of closed curves, namely, circles, convex curves and orthoconvex curves. For circles, the existing results in the literature imply very close upper and lower bounds. The case of convex curves is settled by an easy argument; where as the more non-trivial case of strictly convex curves comes as a consequence of a result on ``grid peeling''~\cite{har}. Next we considered covering by orthoconvex curves. It seems that the case of general orthoconvex curves is difficult. 
\begin{prob}
	Improve the trivial lower and upper bounds, $\lceil (n+1)/4\rceil$ and $\lfloor n/2\rfloor$ respectively, of the minimum number of orthoconvex curves required to cover an $n\times n$ grid.
\end{prob}
So we focused on the simplest orthoconvex curves, namely those with at most one or two inner corners. We showed that at least $2n/5$ (which is achieved for $n=5$) orthoconvex curves with at most one inner corner and $2n/7$ curves with at most two inner corners are required to cover an $n\times n$ grid. 
\begin{prob}
	Is the bound $2n/7$ of Theorem~\ref{thm:2corner} tight? Can one improve it?
\end{prob}
We leave it as an open problem to figure out what happens when there are more inner corners.
\begin{prob}
	What can be said about the minimum number of orthoconvex curves with at most $k$ inner corners, where $k\geq 3$, required to cover an $n\times n$ grid?
\end{prob}
\remove{ Next we looked at covering by non-congruent curves, where we were able to apply results and ideas from incidence geometry. Here we made a conjectural statement on covering by circles of different radii, which came as a consequence of the conjectured bounds on the number of point-circle incidences. Finally, we considered covering by ``small'' curves, i.e. curves with a constant number of grid points. A key ingredient that was used here was the existence of a tiling of $\mathbb{Z}^2$ by translates of these curves. An interesting question that could be asked here is: For which small curves do such tiling exists? Note that, existence of such tiling will give us the minimum number of such curves required to cover an $n\times n$ grid.} Lastly, we mention that in this article we only considered 1-fold covering where every grid point was covered at least once. But, in general, we could ask analogous questions for $r$-fold covering (i.e. we require that every point is covered at least $r$ times) for $r\geq 2$. We feel that answering such questions will be equally interesting and challenging.

\remove{An extended version of this work, with some additional results, has been uploaded to arXiv~\cite{arxiv}, which the interested readers are welcome to have a look.}


\bibliographystyle{abbrvnat}
\bibliography{geom-arbor}
%


\label{sec:biblio}

\end{document}